\documentclass[letterpaper, 10 pt, conference]{ieeeconf}  

\IEEEoverridecommandlockouts                              

\usepackage{siunitx}

\usepackage[hidelinks, bookmarks=false]{hyperref}

\usepackage{graphics} 
\usepackage{epsfig} 
\usepackage{mathptmx} 
\usepackage{times} 
\usepackage{amsmath} 
\usepackage{amssymb}  
\usepackage{cite}

\usepackage{amsthm}  
\newtheorem{theorem}{Theorem}

\let\labelindent\relax
\usepackage{enumitem}

\newtheorem{definition}{Definition}

\newcommand{\Fix}{\operatorname{Fix}}

\renewcommand{\authorrefmark}[1]{\ensuremath{^{#1}}}

\title{\LARGE \bf
    Neuromodulation supports robust rhythmic pattern transitions \\in degenerate central pattern generators with fixed connectivity
}

\author{Arthur Fyon\authorrefmark{1}, Alessio Franci\authorrefmark{1,2}, Pierre Sacré\authorrefmark{1}, and Guillaume Drion\authorrefmark{1}
\thanks{This work was supported by the Belgian Government through the Federal Public Service Policy and Support (BOSA) grant NEMODEI. A.\ Fyon is a Postdoctoral Researcher of the Fonds de la Recherche Scientifique--FNRS. \authorrefmark{1}A.\ Fyon, A.\ Franci, P.\ Sacré, and G.\ Drion are with the Department of Electrical Engineering and Computer Science, University of Liege, 4000 Liege, Belgium. \authorrefmark{2}A. Franci is with the WEL-T Department, WEL Research Institute, 1300 Wavre, Belgium. \{afyon, afranci, p.sacre, gdrion\}@uliege.be.}
}
\begin{document}

\maketitle
\thispagestyle{empty}
\pagestyle{empty}

\begin{abstract}
Many essential biological functions, such as breathing and locomotion, rely on the coordination of robust and adaptable rhythmic patterns, governed by specific network architectures known as connectomes. Rhythmic adaptation is often linked to slow structural modifications of the connectome through synaptic plasticity, but such mechanisms are too slow to support rapid, localized rhythmic transitions. Here, we propose a neuromodulation-based control architecture for dynamically reconfiguring rhythmic activity in networks with fixed connectivity. The key control challenge is to achieve reliable rhythm switching despite neuronal \emph{degeneracy}, a form of structured variability where widely different parameter combinations produce similar functional output. Using equivariant bifurcation theory, we derive necessary symmetry conditions on the neuromodulatory projection topology for the existence of target gaits. We then show that an adaptive neuromodulation controller, operating in a low-dimensional feedback gain space, robustly enforces gait transitions in conductance-based neuron models despite large parametric variability. The framework is validated in simulation on a quadrupedal gait control problem, demonstrating reliable gallop-to-trot transitions across 200~degenerate networks with up to fivefold conductance variability.
\end{abstract}

\section{Introduction}
Generating and switching between rhythmic patterns is a fundamental control problem that arises in both biological and engineered systems. In biology, Central Pattern Generators (CPGs) are self-organized neuronal circuits capable of generating rhythmic motor output in the absence of rhythmic input~\cite{selverston2009}, coordinating vital functions such as locomotion, breathing, and digestion~\cite{marder2001}. A defining feature of biological CPGs is their ability to rapidly switch between rhythmic patterns, for instance, when a quadruped transitions from trot to gallop~\cite{dutta2019}. 
These transitions are typically driven by neuromodulatory inputs from higher brain regions, which dynamically reconfigure network activity by targeting ion channels and synapses~\cite{calabrese1996,bargmann2013}.

In most artificial CPG implementations, rhythm transitions are achieved by altering the network connectivity (the connectome), effectively selecting a different dynamical system for each desired pattern~\cite{barron2010, bay2007, yu2013}. While effective in classical implementations, this strategy faces three fundamental limitations in the context of neuromorphic hardware and biologically plausible control. 

First, connectome modification maps onto synaptic plasticity, a process that operates on timescales of minutes to hours~\cite{morrison2008}, far too slow for the rapid transitions observed in animal locomotion. 

Second, this approach is not robust to \emph{neuronal degeneracy}---the property that widely different parameter configurations can produce equivalent functional output---which is ubiquitous in both biological neural circuits~\cite{goaillard2021, prinz2017} and neuromorphic hardware~\cite{borkar2006} (formally defined in Section~\ref{sec:problem}). Because degenerate neurons respond unpredictably to identical perturbations~\cite{marder2014, fyon2024, fyonnmodhomeo2024, brandoit2025fast}, open-loop connectome switching is an unreliable control strategy. 

Third, these connectome-based approaches rely predominantly on oscillator models, which describe neurons as coupled nonlinear oscillators~\cite{ijspeert2008} and generate different rhythms through connectivity changes rather than local mechanisms such as neuromodulation~\cite{haddad2018}. Oscillator models are widely used in CPG-based robotic systems due to their low computational complexity~\cite{bay2007, yu2013}. However, experimental evidence indicates that biological CPGs achieve adaptive motor control by directly modulating neuronal properties rather than rewiring~\cite{bucher2015}. Capturing this mechanism requires shifting to single-neuron dynamics, specifically conductance-based models following the Hodgkin-Huxley formalism~\cite{hodgkin1952}, which represent neurons as nonlinear RC circuits capturing the ionic channel composition of the neuronal membrane~\cite{mccrea2007}.

These three limitations raise a core question: \emph{how can a network reliably switch between distinct periodic orbits without modifying its connectivity?}

To address this question, we introduce a neuromodulated CPG architecture that achieves adaptive gait control on a fixed connectome, compatible with neuromorphic implementation~\cite{indiveri2011, rubino2019}. Unlike centralized control approaches that rely on high-level optimization or dense sensing~\cite{hoffmanadvancements2022, kotha2024}, our architecture operates through local, state-dependent modulation of neuronal properties. Our contributions are threefold:
\begin{itemize}
    \item We propose a modular CPG architecture with neuromodulatory projections that enables adaptive gait generation for quadruped locomotion through control of neuronal properties rather than connectome modification;
    \item Using equivariant bifurcation theory~\cite{golubitsky2002}, building on the coupled cell framework of Golubitsky, Stewart, and Collins~\cite{golubitsky1998, golubitsky2002, collins1993}, we derive necessary symmetry conditions on the neuromodulatory projection topology for the existence of gaits, revealing that an asymmetric single-neuron projection cannot support the desired periodic orbits, while a symmetric half-center arrangement satisfies all existence conditions;
    \item Through numerical simulations with conductance-based neuron models, we demonstrate reliable gait transitions across 200 degenerate networks with up to fivefold conductance variability, whereas the asymmetric topology produces no stable gait.
\end{itemize}

The remainder of the paper is organized as follows. Section~\ref{sec:problem} formalizes the control problem and defines neuronal degeneracy and robust rhythm switching. Section~\ref{sec:arch} presents the proposed neuromodulated CPG architecture. Section~\ref{sec:sym} derives symmetry conditions for gait existence using equivariant bifurcation theory. Section~\ref{sec:implementation} details the conductance-based neuronal implementation. Section~\ref{sec:design} summarizes the resulting design principles. Section~\ref{sec:casestudy} demonstrates robust gait transitions in simulation.

\section{Problem formulation}\label{sec:problem}

\subsection{Formal setup and the challenge of degeneracy}

We consider a network of $n$ neurons coupled through a fixed connectivity graph (the \emph{connectome}). Each neuron~$i$ is a $k$-dimensional dynamical system with state $x_i \in \mathbb{R}^k$ and the full network state is $x = (x_1,\ldots,x_n) \in \mathbb{R}^{nk}$. The dynamics are given by 
\begin{equation}\label{eq:network_gen}
\dot{x}_i = f(x_i, \theta_i) + \sum_{j \in \mathcal{N}(i)} h(x_j,x_i),  \qquad i = 1, \dots, n,
\end{equation} 
where $f$ describes the intrinsic dynamics of neuron $i$, parameterized by $\theta_i \in \Theta \subseteq \mathbb{R}^p$, and $h$ is a coupling function whose topology is determined by the neighborhood structure $\mathcal{N}(i)$ of the connectome. The connectivity $\mathcal{N}(i)$ and the coupling function $h$ remain fixed throughout; only (a subset of) the intrinsic parameters $\theta_i$ are subject to control. The specific form of $f$ and $h$ is left unspecified at this stage; a conductance-based implementation is detailed in Section~\ref{sec:implementation}.

A central challenge in controlling such networks is that the intrinsic dynamics $f$ exhibit \emph{neuronal degeneracy}: the mapping from parameters to functional output is many-to-one.

\begin{definition}[Neuronal degeneracy]
\label{def:degeneracy}
The intrinsic dynamics $\dot{x}_i = f(x_i, \theta_i)$ exhibit \emph{degeneracy} with respect to an equivalence relation $\sim$ on single-neuron trajectories (\textit{e.g.}, similar oscillation frequency, firing pattern, or burst duty cycle) if there exist $\theta, \theta' \in \Theta$ with $\|\theta - \theta'\| \gg 0$ such that $x(\cdot; \theta) \sim x(\cdot; \theta')$.
\end{definition}

Degeneracy is ubiquitous in biological neural circuits~\cite{goaillard2021, prinz2017} and has a direct analog in neuromorphic hardware through transistor-level manufacturing variability~\cite{borkar2006}. While degeneracy contributes to robustness of baseline function~\cite{goaillard2021}, it poses a fundamental challenge for network control: because neurons with very different intrinsic parameters produce equivalent isolated behavior, their responses to identical perturbations (such as the parameter shifts imposed by connectome switching) are highly variable and unpredictable~\cite{marder2014, fyon2024, fyonnmodhomeo2024, brandoit2025fast}. This makes open-loop connectome modification an unreliable strategy for rhythm transitions in degenerate networks.

\subsection{Control objective}
We seek a closed-loop mechanism that can switch between target gaits on a fixed connectome despite parametric variability across neurons.

\begin{definition}[Robust rhythm switching]\label{def:robust}
Two different network trajectories $x(\cdot)$ and $\tilde{x}(\cdot)$ are \emph{gait-equivalent}, written $x \simeq \tilde{x}$, if they exhibit the same inter-limb phase pattern: the same limbs are synchronized and the same limbs are in antiphase. Consider a family of networks of the form~\eqref{eq:network_gen} with fixed coupling $(h, \mathcal{N})$ and intrinsic parameters $\theta \in \Theta_\mathrm{deg} \subset \mathbb{R}^p$, where $\Theta_\mathrm{deg}$ is a degenerate population (Definition~\ref{def:degeneracy}). A control mechanism achieves \emph{robust rhythm switching} if, for all $\theta \in \Theta_\mathrm{deg}$, the closed-loop system converges to a periodic orbit that is gait-equivalent ($\simeq$) to the target gait within a bounded transient time.
\end{definition}

The central question addressed in this paper is therefore: given a fixed connectome of neurons with degenerate parameters, design a control architecture that achieves robust rhythm switching (Definition~\ref{def:robust}) without modifying the network connectivity.

\subsection{Case study: quadrupedal gaits}

\begin{figure*}[!ht]
  \vspace{0.2cm}
  \centering
  \includegraphics[width=0.7\linewidth]{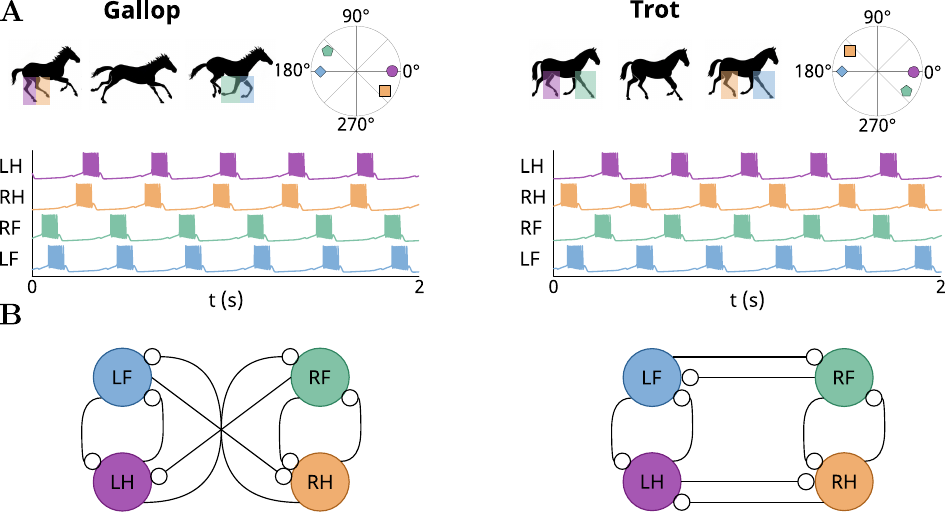}
  \caption{\textbf{Gait: gallop and trot.} \textbf{A.} Trot and gallop rhythms with associated phase patterns (top) and the corresponding bursting behavior of motor neurons controlling each leg (bottom). \textbf{B.} State-of-the-art artificial neuronal structure with distinct connectomes for each rhythm. White dots represent inhibitory synapses. Labels: L for left, R for right, F for front, and H for hind.}
  \label{fig:1}
\end{figure*}

We demonstrate our approach on a simplified quadruped model that switches between two rhythmic patterns: gallop and trot (Fig.~\ref{fig:1}A). During each gait cycle, muscles contract during the active phase of the corresponding motor neuron bursts, ensuring coordinated limb movement.

These rhythms are intuitive: in the trot gait, diagonally opposed legs alternate, while in the gallop gait, the front or hind legs move in unison. The coordination between motor neuron bursting and limb activation can be observed in the phase patterns (Fig.~\ref{fig:1}A), highlighting the structured timing of muscle activation. 

In artificial neuronal systems, transitions between these rhythms are typically achieved by modifying the network connectome (Fig.~\ref{fig:1}B)~\cite{barron2010}. Our goal is to achieve the same transitions through neuromodulatory control alone.

\section{Proposed architecture}\label{sec:arch}

\begin{figure}[t!]
  \centering
  \includegraphics[scale=0.37]{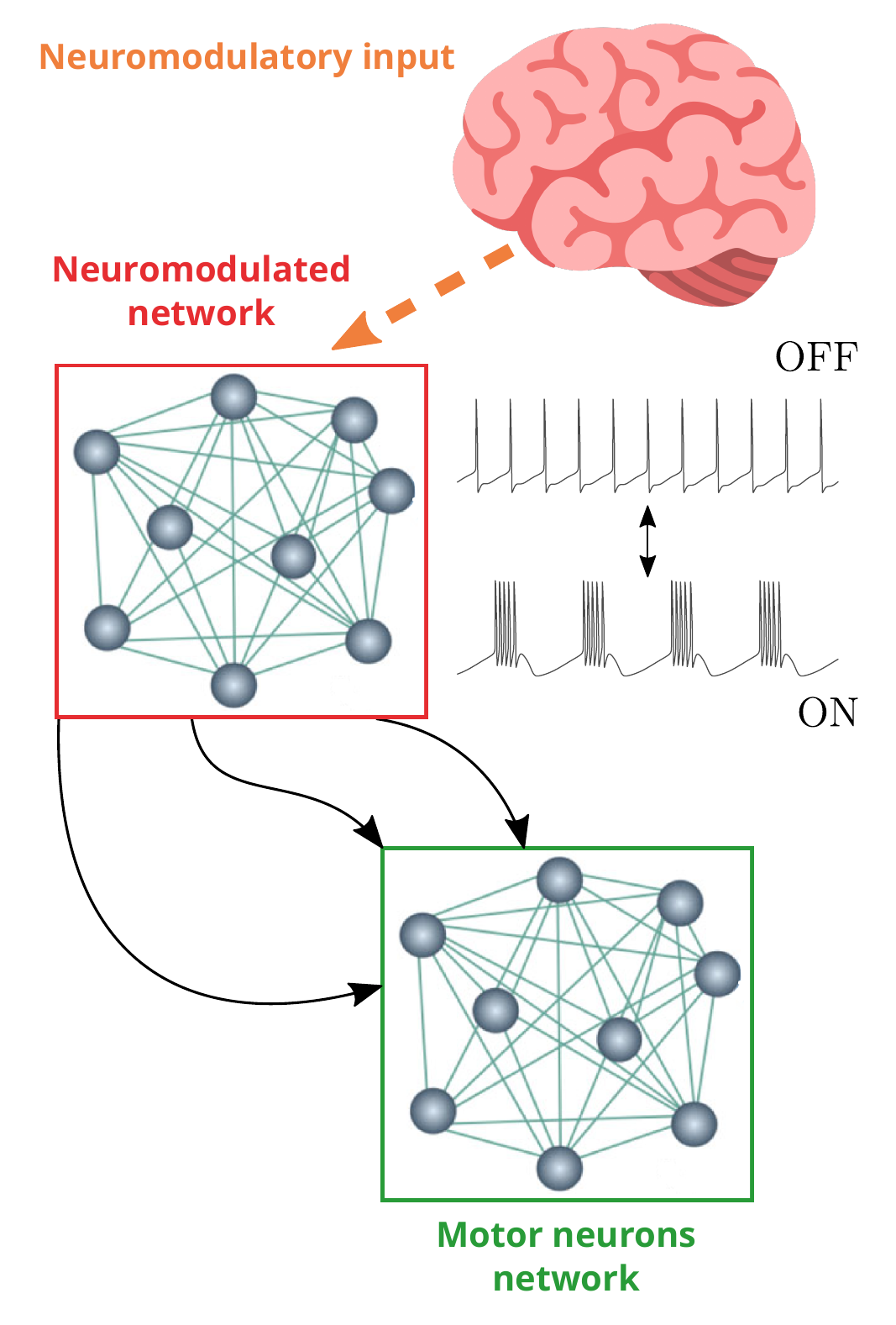}
  \caption{\textbf{The neuromodulated central pattern generator architecture.} The top-down organization of the neuromodulated CPG architecture, showing the interaction between the neuromodulatory input, the neuromodulated network, and the motor neuron network.}
  \label{fig:2}
\end{figure}

The proposed neuromodulated CPG architecture consists of three interconnected layers (Fig.~\ref{fig:2}), mirroring the hierarchical organization of biological locomotion control~\cite{dutta2019, bucher2015, engel2008}. The switching principle is model-agnostic: it applies to any neuron model that supports transitions between tonic firing and bursting through parameter modulation. The specific conductance-based implementation used for validation is detailed in Section~\ref{sec:implementation}.

Building on this organization, the three layers are:
\begin{itemize}
    \item \textbf{A neuromodulatory input}, which projects to the neuromodulated network. The neuromodulatory controller~\cite{fyon2023} provides state-dependent adaptation of neuronal properties, enabling neurons to transition between tonic firing and bursting while maintaining robustness to degeneracy;
    
    \item \textbf{A neuromodulated network}, composed of bursting neurons with inhibitory synapses organized as half-center oscillators, a fundamental CPG motif found in both biological and artificial systems~\cite{calabrese1995, sepulchre2022}. This network projects inhibitory synapses onto the motor neurons, and neuromodulation alters the strength of these projections to enable rhythmic pattern transitions;
    
    \item \textbf{A motor neuron network}, consisting of bursting neurons that drive muscle activation and include mutual inhibitory synapses.
\end{itemize}

Rhythmic transitions occur through the modulation of specific neurons within the neuromodulated network, which switch between bursting and tonic firing states (silent activity is excluded as it is incompatible with homeostatic tuning rules~\cite{oleary2014}). Importantly, the network does not rely on synaptic plasticity to achieve these transitions. The rhythmic patterns are embedded in the projections from the neuromodulated network to the motor neuron network: when a given rhythm is required, only the corresponding neurons remain in bursting mode while the others switch to tonic firing, so that the dominant inhibitory projections from the bursting neurons enforce the desired pattern on the motor layer.

\section{Symmetry-based existence of gaits}\label{sec:sym}
The architecture proposed in Section~\ref{sec:arch} raises a fundamental design question: \emph{which topologies of neuromodulatory projections can support the desired gaits?} We address this question using equivariant bifurcation theory, building directly on the coupled cell framework~\cite{golubitsky2002, golubitsky1998, collins1993} and applying it to our specific architecture.

\begin{figure*}[!ht]
  \vspace{0.2cm}
  \centering
  \includegraphics[width=0.85\linewidth]{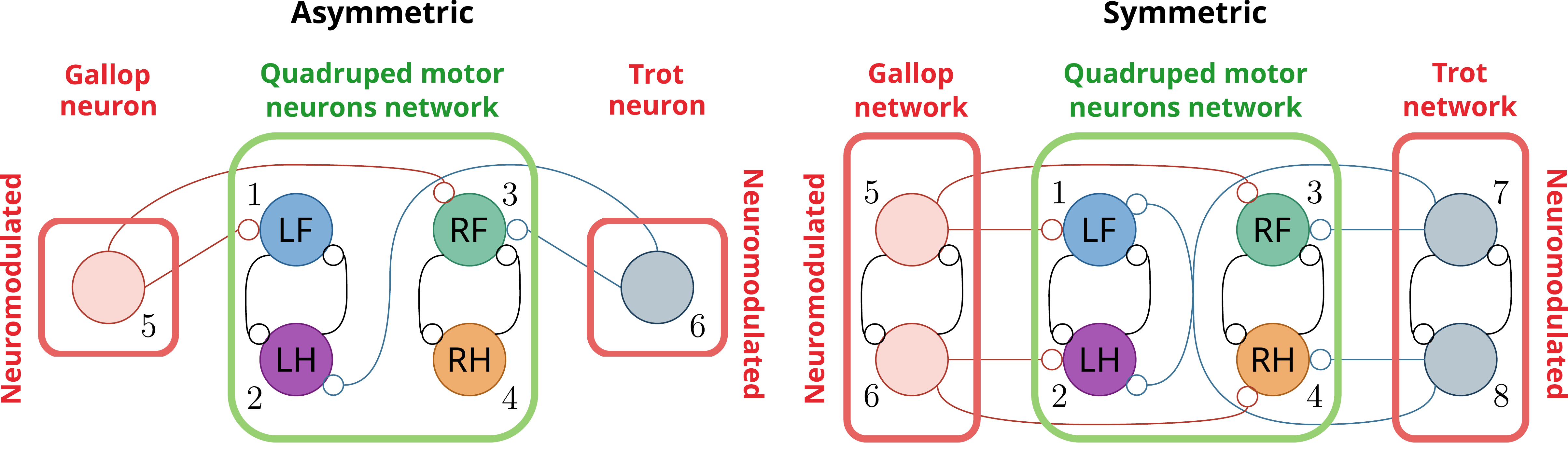}
  \caption{\textbf{Different studied topologies.} Minimal asymmetric architecture including four motor neurons and two neurons in the neuromodulated network, one for each rhythm in the toy gait control example (left). Minimal symmetric architecture including four motor neurons and four neurons in the neuromodulated network, two for each rhythm in the toy gait control example (right).}
  \label{fig:topo}
\end{figure*}

The motor network architecture follows a natural and well-established organization, where distinct motor neurons encode each rhythm. We analyze two minimal neuromodulatory topologies (Fig.~\ref{fig:topo}): an asymmetric architecture with a single neuron encoding each rhythm (left), and a symmetric architecture with two neurons per rhythm arranged in half-center oscillators (right). For the symmetry analysis, we focus on the gallop rhythm: gallop-coding neurons remain in bursting mode, while trot-coding neurons switch to tonic firing. Their inhibitory outputs can be discarded, reducing the left topology to five nodes (1--5) and the right topology to six nodes (1--6).

\textbf{Mathematical framework.} 
Recall from Section \ref{sec:problem} that the network state lies in $\mathbb{R}^{nk}$ and the dynamics respect the coupled cell structure~\eqref{eq:network_gen}. The network dynamics respect the symmetries of the coupling graph: specifically, the system is \emph{$\Gamma$-equivariant}, where $\Gamma$ denotes the automorphism group of the coupling graph defined by $\mathcal{N}$ (the group of all permutations preserving connectivity), meaning that if $\gamma\in\Gamma$ permutes the nodes, the vector field transforms accordingly. We further denote by $K \leq \Gamma$ a spatial isotropy subgroup, by $H \leq \Gamma$ a spatio-temporal symmetry group, and by $\Fix(K)$ the fixed-point subspace of $K$.

Periodic solutions to such symmetric systems inherit specific symmetry properties. A \emph{spatial symmetry} $K$ is a subgroup of $\Gamma$ that leaves the periodic orbit invariant at each instant in time (\textit{e.g.}, left-right synchrony). A \emph{spatio-temporal symmetry} $H$ is a larger subgroup that leaves the orbit invariant up to a time shift (\textit{e.g.}, a half-period shift that swaps front and hind limbs). The relationship between these symmetries and the existence of periodic orbits is characterized by the following fundamental result from equivariant bifurcation theory.

\begin{theorem}[Theorem~3.4 in~\cite{golubitsky2002}]
\label{thm:equivariant_periodic}
Let $\Gamma$ be a finite group acting on $\mathbb{R}^{n}$. There exists a $T$-periodic solution to a $\Gamma$-equivariant system of ODEs on $\mathbb{R}^{n}$ with spatial symmetries $K$ and spatio-temporal symmetries $H$ if and only if:
\begin{enumerate}[label=(\alph*)]
    \item $H/K$ is cyclic;
    \item $K$ is an isotropy subgroup;
    \item $\dim \Fix(K)\geq 2$ (with extra conditions if equality holds);
    \item $H$ fixes a connected component of $\Fix(K)\setminus L_K$, where $L_K=\bigcup_{\gamma\notin K} \Fix(\gamma)\cap\Fix(K)$.
\end{enumerate}
\end{theorem}

The fixed-point subspace of a subgroup $K\leq \Gamma$ is denoted
\[
\Fix(K)=\{x\in\mathbb R^n:\;\gamma x=x\ \text{for all }\gamma\in K\}.
\]

\subsection{Nonexistence of a gallop with a single projection neuron}
Consider a 5-node network with four limb nodes and a single projection neuron $5$ (Fig.~\ref{fig:topo} left). Define the two involutions
\[
a=(1\;2)(3\;4),\qquad b=(1\;3)(2\;4)(5),
\]
so the symmetry group is $\Gamma=\langle a,b\rangle\cong \mathbb Z_2\times\mathbb Z_2$.

For a gallop, the natural spatial isotropy is $K=\langle a\rangle=\{e,a\}$, enforcing left-right synchrony within the front and hind pairs. The fixed-point subspace is
\[
\Fix(K)=\{x\in\mathbb R^{5k}:\;x_1=x_2,\;x_3=x_4\}\cong \mathbb R^{3k},
\]
with coordinates $(u=x_1=x_2, v=x_3=x_4, w=x_5)$.

The elements outside $K$ ($b$ and $c=ab$) impose $u=v$, so
\[
L_K=\{u=v\}\subset\Fix(K).
\]
This is codimension one, hence
\[
\Fix(K)\setminus L_K=\{(u,v,w)\in\mathbb R^{3k}:\ u\neq v\},
\]
which has two connected components ($u>v$ and $u<v$). Any spatio-temporal group $H=\Gamma$ containing $b$ necessarily exchanges these two components rather than preserving one. Thus condition (d) fails, and no gallop solution exists with this topology.

\subsection{Existence of a gallop with two projection neurons}

Now consider the 6-node network with two symmetric projection neurons (Fig.~\ref{fig:topo} right). Define
\[
a=(1\;2)(3\;4),\qquad b=(1\;3)(2\;4)(5\;6),
\]
so the group is $\Gamma=\langle a,b\rangle\cong \mathbb Z_2\times\mathbb Z_2$ with $c=ab=(1\;4)(2\;3)(5\;6)$.

The spatial isotropy for a gallop is again $K=\langle a\rangle=\{e,a\}$, and the fixed-point subspace is
\[
\Fix(K)=\{x\in\mathbb R^{6k}:\;x_1=x_2,\;x_3=x_4\}\cong \mathbb R^{4k},
\]
with coordinates $(u=x_1=x_2, v=x_3=x_4, w=x_5, z=x_6)$.

The extra symmetries $b$ and $c$ impose $u=v$ and $w=z$, so
\[
L_K=\{u=v,\;w=z\},
\]
a codimension-two subspace. Removing it leaves $\Fix(K)\setminus L_K$ connected. The spatio-temporal group $H=\Gamma=\{e,a,b,c\}$ preserves this connected component, since $b$ corresponds to a half-period time shift swapping front and hind. All conditions (a)--(d) of Theorem~\ref{thm:equivariant_periodic} hold, so a gallop solution exists.

\subsection{Extension to trot}

The same reasoning applies to the trot, where the isotropy subgroup enforces diagonal synchrony ($x_1=x_4$, $x_2=x_3$). Neuromodulation selectively activates projection neurons, thereby switching between symmetric topologies that support trot or gallop without requiring plasticity.

\section{Neuronal implementation}\label{sec:implementation}

Having established which projection topologies can support the desired gaits, we now turn to the second design question: \emph{what control strategy can reliably realize gait transitions despite neuronal degeneracy?} 

We instantiate the framework of Section~\ref{sec:problem} using conductance-based neuron models. We focus on biophysical models with the aim of demonstrating robustness at this level first, so that the same control principle can later be translated to simplified neuron models suited for mixed-feedback neuromorphic implementations operating in the nanowatt power range.

\subsection{Conductance-based neuron models}
In single-compartment conductance-based models, each neuron intrinsic dynamics $f(x_i, \theta_i)$ in~\eqref{eq:network_gen} take the form of a voltage-current relationship:
\begin{align}
   \underbrace{C \frac{\mathrm{d}V}{\mathrm{d}t} + g_\mathrm{l}(V-E_\mathrm{l})}_{I_C} = -\mathord{\underbrace{\sum_{\mathrm{ion} \in \mathcal{I}} g_\mathrm{ion}(V,t) (V - E_{\mathrm{ion}})}_{I_\mathrm{int}}} + I_\mathrm{ext}, \label{eq:CBmodel}
\end{align}
where $I_C$ is the membrane current through leakage channels and $I_\mathrm{int}$ is the current through voltage- and time-dependent ion channels (see~\cite{hodgkin1952} for details). Each channel type has a conductance $g_\mathrm{ion}(V,t)$ varying between $0$ and a maximum value $\bar{g}_\mathrm{ion}$ determined by channel density. Importantly, the parameter vector $\bar{g}_\mathrm{ion} \in \mathbb{R}^p$ corresponds to $\theta$ in Definition~\ref{def:degeneracy}: degeneracy manifests as multiple, widely separated configurations of $\bar{g}_\mathrm{ion}$ producing qualitatively equivalent firing patterns.

We employ a stomatogastric (STG) neuron model~\cite{liu1998} with seven ion channel types (Na, CaT, CaS, A, KCa, Kd, H), selected for its biological relevance in rhythmic activity and bursting behavior~\cite{bucher2006}. While conductance-based models are state-of-the-art for rhythmic systems, our architecture is intended to be general. To incorporate synaptic connections, we modify \eqref{eq:CBmodel} with a synaptic current term: $I_C = -I_\mathrm{int} - I_\mathrm{syn} + I_\mathrm{ext}$, where the synaptic current between presynaptic voltage $V_\mathrm{pre}$ and postsynaptic voltage $V_\mathrm{post}$ is:
\begin{equation}
    I_\mathrm{syn}(V_\mathrm{pre}, V_\mathrm{post}, t) = g_\mathrm{syn}(V_\mathrm{pre}, t)(V_\mathrm{post} - E_\mathrm{syn}),
\end{equation}
with $g_\mathrm{syn}$ following the classical Hodgkin-Huxley formalism and gating kinetics modeled as in~\cite{drion2019}. For the inhibitory GABAergic synapses considered here, the synaptic reversal potential is approximately $\SI{-75}{mV}$. This synaptic current instantiates the coupling function $h$ in~\eqref{eq:network_gen}.

\subsection{Neuromodulation as adaptive control}
Neuromodulation tunes the maximum conductances $\bar{g}_\mathrm{ion}$ to transition neurons between tonic firing and bursting~\cite{zagha2014}. State-independent, fixed modulation rules are not robust to degeneracy~\cite{marder2014}; a state-dependent approach is required.

\begin{figure}[t!]
  \vspace{0.2cm}
  \centering
  \includegraphics[scale=0.95]{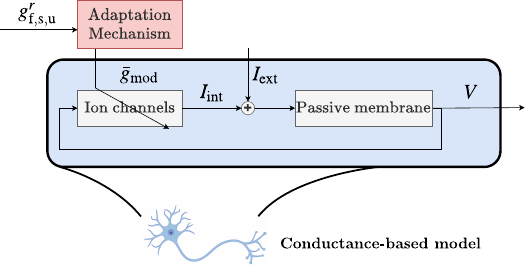}
  \caption{\textbf{The adaptive neuromodulation controller.} High-level block diagram of the adaptive neuromodulation controller. The blue block represents the typical structure of a conductance-based model from a feedback control perspective. The red block groups all biological mechanisms that regulate ion channel expression and function as an adaptive layer on the neuronal controller. This adaptive mechanism takes target neuronal feedback gains, hence target activity, as input and adjusts a subset of ion channels $\bar{g}_\mathrm{mod}$ to achieve the desired activity. Adapted from~\cite{fyon2023}.}
  \label{fig:1b}
\end{figure}

The framework of~\cite{fyon2023} formulates neuromodulation as an adaptive control problem, viewing the conductance-based model as a feedback system (Fig.~\ref{fig:1b}): ion channels act as a neuronal controller and the passive membrane acts as the plant.

To enable tractable control, all ion channels are aggregated into three voltage-dependent neuronal feedback gains, termed Dynamic Input Conductances (DICs), which are separated by timescales~\cite{drion2015dic}. While the raw conductances $\bar{g}_\mathrm{ion}$ do not directly predict neuronal behavior, the feedback gains evaluated at the threshold voltage $V_\mathrm{th}$ provide an interpretable and predictive description of spiking activity. 
The neuronal feedback gains are computed as:
\begin{equation}\label{eq:DICsystem}
\begin{bmatrix}
g_\mathrm{f}(V_\mathrm{th})\\
g_\mathrm{s}(V_\mathrm{th})\\
g_\mathrm{u}(V_\mathrm{th})
\end{bmatrix}= S(V_\mathrm{th}) \cdot \bar{g}_\mathrm{ion}\, ,
\end{equation}
where $S(V_\mathrm{th}) \in \mathbb{R}^{3 \times p}$ is a sensitivity matrix computed as in~\cite{drion2015dic}. Distinct regions of this three-dimensional gain space correspond to tonic firing, bursting, and silence.

Following the principles of Model Reference Adaptive Control (MRAC)~\cite{MRAC}, the neuromodulation controller takes as input target neuronal feedback gains: 
$$g^r_\mathrm{f, s, u} = (g^r_\mathrm{f}(V_\mathrm{th}), g^r_\mathrm{s}(V_\mathrm{th}), g^r_\mathrm{u}(V_\mathrm{th})),$$
which define the desired spiking activity. The target gains are set by the higher-level gait selection logic, specifying whether each neuron should be in bursting or tonic mode. The controller then autonomously adjusts a subset of $\bar{g}_\mathrm{ion}$, $\bar{g}_\mathrm{mod}$, in a state-dependent manner to drive the actual neuronal feedback gains toward the target values (Fig.~\ref{fig:1b}). By operating within an interpretable and tractable space of firing activity, the controller effectively adapts conductances in the complex and high-dimensional space of $\bar{g}_\mathrm{ion}$ by solving \eqref{eq:DICsystem}. Because the controller is closed-loop and state-dependent, it adapts its actions to each neuron current parameter configuration, regardless of where in the degenerate manifold that neuron lies, enabling neuromodulatory effects that are robust to degeneracy.

This formulation has a direct neuromorphic analog: in silicon implementations, bias currents (analogous to $\bar{g}_\mathrm{mod}$) compensate for manufacturing mismatch (analogous to uncontrollable $\bar{g}_\mathrm{ion} \setminus \bar{g}_\mathrm{mod}$), providing a principled decomposition into controllable and uncontrollable parameters that justifies modulating only a subset of conductances.

Although the neuromodulation controller was introduced in prior work as a single-neuron mechanism~\cite{fyon2023}, the present study demonstrates its integration into a structured network architecture. This application to network-level pattern reconfiguration, robustly enabling rhythm transitions in the presence of degeneracy, constitutes a key extension beyond the scope of~\cite{fyon2023}.

\section{Design principles}\label{sec:design}

The combined symmetry analysis and control formulation yield concrete design requirements for our architecture:

\begin{enumerate}[label=(\roman*)]
    \item \textbf{Projection symmetry:} The neuromodulated layer must have symmetric projections (half-center oscillators) encoding each rhythm, as established by the symmetry analysis (Section~\ref{sec:sym}). Asymmetric single-neuron projections do not support the required gaits.
    
    \item \textbf{Rich neuronal dynamics:} 
    Individual neurons must support multiple qualitatively distinct activity regimes (at minimum tonic firing and bursting) that can be accessed through parameter modulation. This is satisfied by any neuron model exhibiting transitions between tonic and bursting regimes as a function of tunable parameters~\cite{drion2015}. In degenerate populations, the transition boundary in parameter space varies across neurons. Combined with \eqref{eq:DICsystem}, these rich dynamics ensure that the controller can steer neuronal feedback gains between the tonic and bursting regions.
    
    \item \textbf{Coupling strength:} Synaptic connections use inhibitory coupling. The maximum synaptic conductance $\bar{g}_\mathrm{syn}$ must be large enough for the projections from bursting neurons in the neuromodulated layer to entrain the motor layer, but not so large as to suppress all motor neuron activity. In our simulations, $\bar{g}_\mathrm{syn} = \SI{0.8}{mS/cm^2}$ satisfies this constraint across the full degenerate population.
\end{enumerate}

These requirements are structural rather than fine-tuned: they constrain the topology and the qualitative capabilities of the neuron model, not specific parameter values. This makes the framework compatible with a broad class of spiking neuron implementations, including mixed-feedback neuromorphic circuits~\cite{indiveri2011, rubino2019}, where target gains map to bias currents set by an upstream control circuit.

\section{Robust gait control (case study)}\label{sec:casestudy}

\begin{figure*}[!t]
\vspace{0.2cm}
\centering\includegraphics[width=0.85\linewidth]{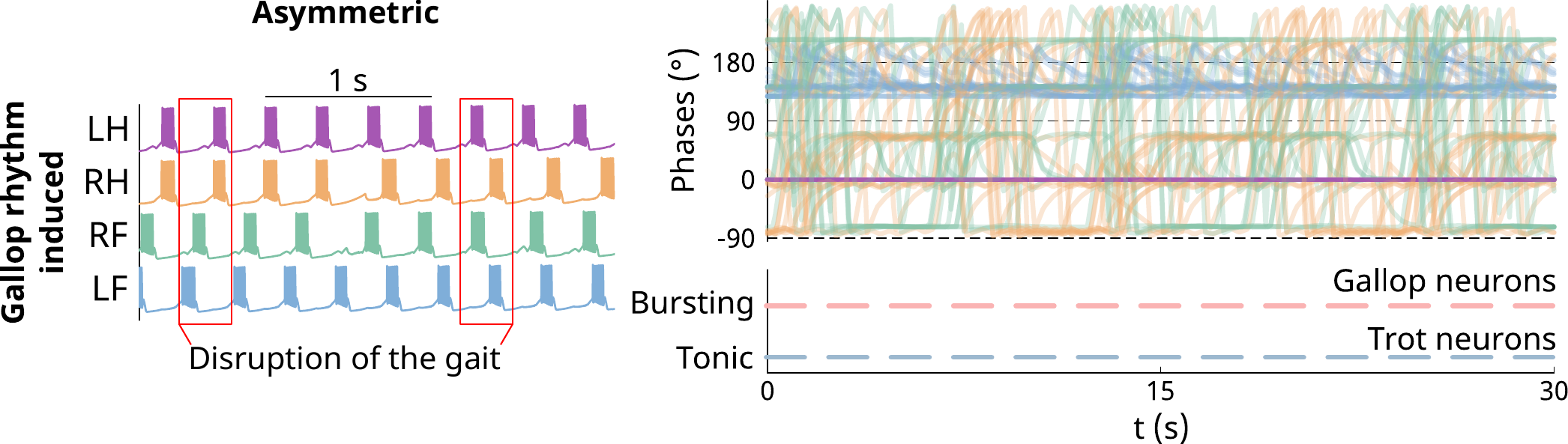}
        \caption{\textbf{Asymmetric gait control leads to gait disruption.} Voltage traces of the four motor neurons during gallop rhythm, without switching, for the asymmetric architecture (left). Phase patterns during the gallop rhythm of the motor neurons for different initial conditions with $V(t=0)$ drawn from a uniform distribution ranging from $\SI{-65}{mV}$ to $\SI{-55}{mV}$ (right).}
        \label{fig:4}
\end{figure*}

\begin{figure*}[!t]
\centering
\includegraphics[width=0.85\linewidth]{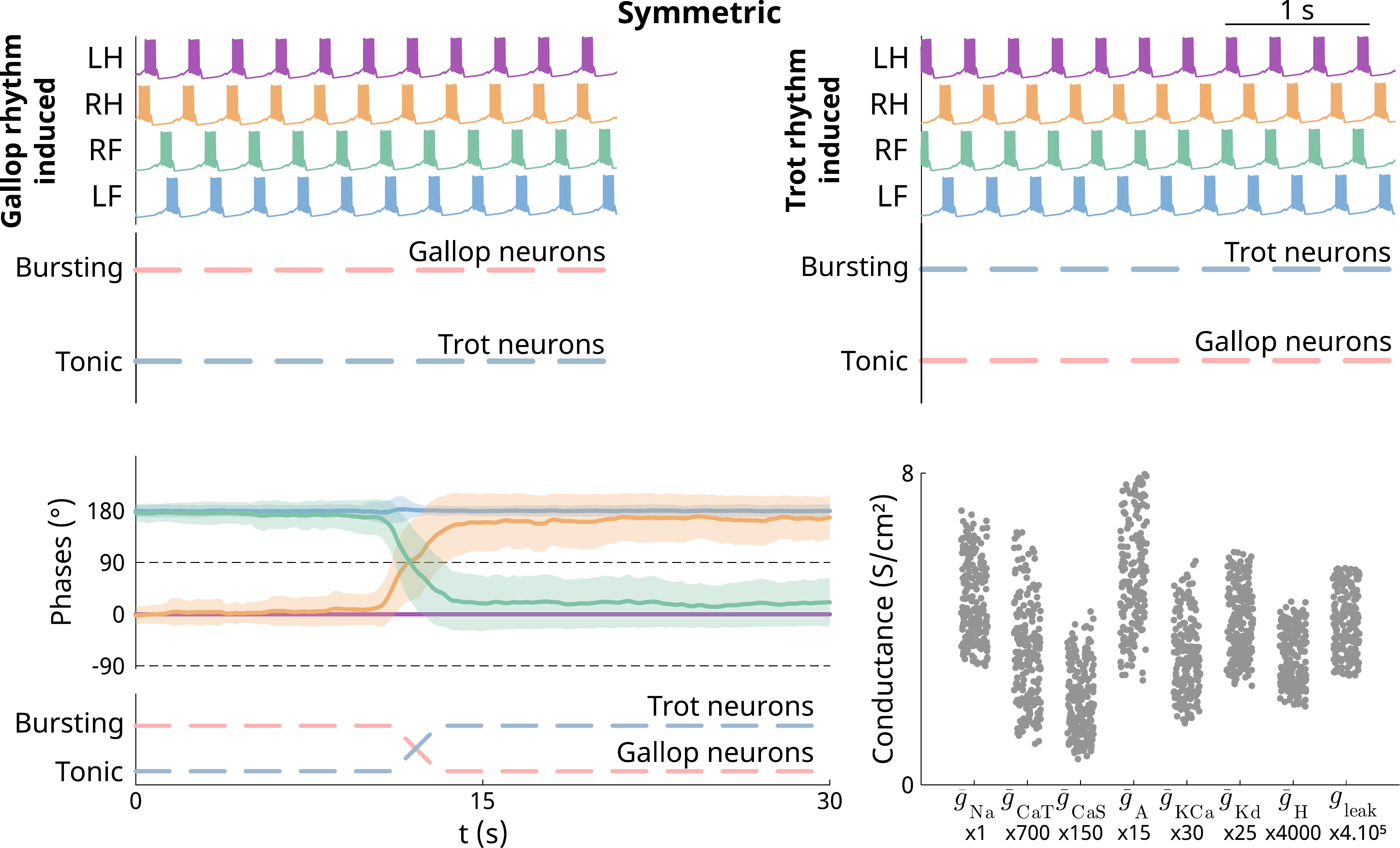}
        \caption{\textbf{Symmetric gait control with variable conductances.} Voltage traces of the four motor neurons during gallop and trot rhythms, without switching, for the symmetric architecture (top). Phase patterns of the motor neurons undergoing the transition from gallop to trot with variable conductances (bottom left). The solid line represents the mean across the population of networks, while the shaded areas indicate the min-max range. Distributions of ionic conductances across the neuronal population far exceeding biologically plausible limits (bottom right).}
        \label{fig:5}
\end{figure*}

Rhythmic output is classified as ``gallop'' or ``trot'' based on the inter-limb phase relationships between bursts of the four motor neurons (Figs.~\ref{fig:4}--\ref{fig:5}). We now verify the theoretical predictions of Sections~\ref{sec:sym}--\ref{sec:design} in simulation.

\subsection{Asymmetric gait control leads to gait disruption}
As demonstrated in Section~\ref{sec:sym}, the asymmetric architecture (Fig.~\ref{fig:topo} left) cannot support stable gallop solutions due to the absence of required symmetry properties. Simulations confirm this theoretical prediction: when the gallop neurons are activated, the network fails to sustain any regular gait pattern. Instead, continuous phase switching occurs, preventing convergence to a periodic solution (Fig.~\ref{fig:4} left). 

This instability persists across 20 simulations with varying initial conditions (Fig.~\ref{fig:4} right), and identical observations hold for the trot rhythm. The symmetry conditions of Section~\ref{sec:sym} are therefore necessary for stable gait generation.

\subsection{Symmetric gait control under degeneracy}
The symmetric architecture (Fig.~\ref{fig:topo} right) reliably generates the expected motor rhythm when either gallop or trot neuromodulated neurons are activated. This projection symmetry ensures stable and sustained rhythmic patterns, consistent with the theoretical predictions of Section~\ref{sec:sym} (Fig.~\ref{fig:5} top). Critically, transitions between rhythms can be achieved without any connectome modifications.

To demonstrate robustness in the sense of Definition~\ref{def:robust}, we generated 200 networks with varying conductance values using the algorithm from~\cite{fyon2024}. Despite having identical connectome structure and synaptic strengths, these networks exhibit up to fivefold variations in intrinsic conductances across neurons. Fig.~\ref{fig:5} (bottom) shows phase patterns during gallop-to-trot transitions, triggered by switching gallop neurons from bursting to tonic spiking while simultaneously activating trot neurons. The symmetric architecture maintains stable rhythmic transitions across all 200 networks, with the qualitative phase pattern (in the sense of gait equivalence $\simeq$, Definition~\ref{def:robust}) preserved across the entire degenerate population.

This modular architecture extends naturally to larger CPGs or more complex locomotion patterns, as each neuromodulated neuron operates independently using local feedback gains. In our framework, muscle activation aligns with motor neuron burst phases, while silent phases correspond to inactivity. 


\section{Conclusions and future work}
Symmetric neuromodulatory projections are necessary for gait generation on a fixed connectome, and adaptive neuromodulation is sufficient for reliable switching despite degeneracy, even under fivefold conductance variability far exceeding biological limits. This result separates the algebraic question (\emph{which} topologies support which gaits, addressed by equivariant bifurcation theory) from the control question (\emph{how} to robustly realize them, addressed by adaptive neuromodulation), and shows that both can be resolved on a single fixed connectome without synaptic plasticity.

Because the controller acts at the neuron level using local feedback rules, this approach is inherently scalable: larger networks can be assembled by adding modules without changing the core control architecture. This method could pave the way for a new generation of bio-inspired complex locomotion control in neuromorphic systems.
                           
Future work will explore how burst timing modulates actuation parameters such as gait speed and force.
This work focuses solely on biophysical models. The underlying theories, equivariant bifurcation theory and DICs, apply to a broad class of models, and future work will generalize these results to mixed feedback systems and CMOS neuromorphic neurons. In particular, the neuromodulation controller does not require full simulation of biophysical models, but can be mapped onto analog control loops that regulate a few key conductance parameters, making it suitable for efficient neuromorphic implementation. Building on this modular CPG approach, complex neuromorphic locomotion systems can then be designed, ultimately leading to validation on a real quadruped robot. Energy profiling and control simplification will be investigated during this process to ensure compatibility with ultra-low-power hardware.


\section*{Acknowledgment}

The Claude Opus 4.6 chatbot was used to improve the syntax and grammar of several paragraphs in this manuscript. The authors subsequently reviewed and edited the text as needed and take full responsibility for the published content.


\section*{Code availability}
All code and data can be found on the first author's GitHub: \href{https://github.com/arthur-fyon/ECC_2026}{\textbf{https://github.com/arthur-fyon/ECC\_2026}}.

\addtolength{\textheight}{-4.9cm}   
                                  

\end{document}